\documentclass{amsart}
\usepackage{amsmath}
\usepackage{amsthm}
\usepackage{latexsym}
\usepackage{amssymb}
\usepackage{xspace}
\usepackage{amscd}
\usepackage[latin1]{inputenc}
\usepackage[dvips]{graphicx}

\date{}

\newcommand{\Gal}{{\rm Gal}}

\newcommand{\inv}{^{-1}}

\newcommand{\Z}{{\mathbb Z}}
\newcommand{\Q}{{\mathbb Q}}
\newcommand{\R}{{\mathbb R}}
\newcommand{\C}{{\mathbb C}}
\newcommand{\HQ}{{\mathbb H}}

\newcommand{\SL}{{\rm SL}}

\newcommand{\matriz}[1]{\begin{array} #1 \end{array}}

\newcommand{\GEN}[1]{\langle #1 \rangle}
\newcommand{\gorro}[1]{\widehat{#1}}
\newcommand{\quat}[2]{\left( \frac{#1}{#2} \right)}
\newcommand{\U}{{\mathcal U}}
\newcommand{\V}{{\mathcal V}}
\newcommand{\T}{{\mathcal T}}
\newcommand{\W}{{\mathcal W}}
\newcommand{\s}{{\mathcal S}}
\newcommand{\CC}{{\mathcal C}}

\newtheorem{theorem}{Theorem}[section]
\newtheorem{lemma}[theorem]{Lemma}
\newtheorem{proposition}[theorem]{Proposition}

\newtheorem{corollary}[theorem]{Corollary}

\theoremstyle{remark}
\theoremstyle{remark}
\theoremstyle{remark}
\theoremstyle{remark}
\theoremstyle{remark}

\begin{document}

\title{Group algebras of Kleinian type and groups of units}

\author{Gabriela Olteanu}

\address{Department of Mathematics and Computer Science,
North University of Baia Mare, Victoriei 76, 430072 Baia Mare, Romania.
\tt{olteanu@math.ubbcluj.ro}} \curraddr{Department of Mathematics, University of Murcia,
30100 Murcia, Spain. \tt{golteanu@um.es}}

\author{\'{A}ngel del R\'{\i}o}
\address{Department of Mathematics, University of Murcia, 30100 Murcia, Spain. \tt{adelrio@um.es}.
http://www.um.es/adelrio}

\thanks{Research supported by M.E.C. of Romania (CEEX-ET 47/2006),
D.G.I. of Spain and Fundaci\'{o}n S\'{e}neca of Murcia.}



\begin{abstract}
The algebras of Kleinian type are finite dimensional semisimple rational algebras $A$ such that the group of units of
an order in $A$ is commensurable with a direct product of Kleinian groups. We classify the Schur algebras of Kleinian
type and the group algebras of Kleinian type. As an application, we characterize the group rings $RG$, with $R$ an
order in a number field and $G$ a finite group, such that $RG^*$ is virtually a direct product of free-by-free groups.
\end{abstract}

\maketitle

The study of Kleinian groups goes back to the works of Poincar\'{e} \cite{Poi} and
Bianchi \cite{Bia} and it has been an active field of research ever since. In the last
decades, it is strongly related to the Geometrization Program of Thurston for the
classification of 3-manifolds \cite{EGM,MR,Mas,Thu}. The use of the methods of Kleinian
groups to the study of the groups of units of group rings was started in \cite{PRR} and
led to the notion of algebras of Kleinian type and finite groups of Kleinian type.

Let $K$ be a number field, $A$ a central simple $K$-algebra and $R$ a $\Z$-order in $A$.
Let $R^1$ denote the group of units of $R$ of reduced norm $1$. Every embedding
$\sigma:K\rightarrow \C$ induces an embedding $\overline{\sigma}:A\rightarrow M_d(\C)$,
where $d$ is the degree of $A$. Furthermore, $\overline{\sigma}(R^1)\subseteq \SL_d(\C)$.
We say that $A$ is of {\it Kleinian type} if either $A=K$ or $A$ is a quaternion algebra
over $K$ and $\overline{\sigma}(R^1)$ is a discrete subgroup of $\SL_2(\C)$ for some
embedding $\sigma$ of $K$ in $\C$. More generally, an {\it algebra of Kleinian type}
\cite{PRR} is by definition a direct sum of simple algebras of Kleinian type. A {\it
finite group} $G$ is of {\it Kleinian type} if and only if the rational group algebra $\Q
G$ is of Kleinian type.

The finite groups of Kleinian type have been classified in \cite{JPRRZ} where it has been
also proved that a finite group $G$ is of Kleinian type if and only if the group of units
$\Z G^*$ of its integral group ring $\Z G$ is commensurable with a direct product of
free-by-free groups. Recently, Alan Reid asked us in a private communication about the
consequences of replacing the ring of rational integers by another ring of integers. This
leads to the following two problems:

\begin{quote}
{\it Problem 1}. Classify the group algebras of Kleinian type of finite groups over number
fields. \\
{\it Problem 2}. Given a group algebra of Kleinian type $KG$, describe the structure of
the group of units of the group ring $RG$ for $R$ an order in $K$.
\end{quote}

The simple factors of $KG$ are Schur algebras over its center. So, in order to solve
Problem 1, it is natural to start classifying the Schur algebras of Kleinian type. This
is obtained in Section 1. Using this classification and that of finite groups of Kleinian
type given in \cite{JPRRZ} we obtain the classification of the group algebras of Kleinian
type in Section 2. In Section 3 we obtain a partial solution for Problem 2.

\section{Schur algebras}

All throughout $K$ is a number field. We refer to the field homomorphisms $K\rightarrow \R$ as real embeddings of $K$.
By abuse of notation, a pair of complex embeddings of $K$ is, by definition, a pair of conjugate field homomorphisms
$K\rightarrow \C$ whose images are not embedded in $\R$. Recall that the infinite places of $K$ correspond to the real
embeddings of $K$ and the pairs of complex embeddings of $K$. Furthermore, the finite places of $K$ correspond to the
prime ideals of the ring of integers of $K$.

A {\it Schur algebra} over $K$ is a central simple $K$-algebra $A$ which is generated
over $K$ by a finite subgroup of the group of units $A^*$ of $A$. Equivalently, a Schur
algebra over $K$ is a simple factor, with center $K$, of a group algebra of a finite
group.

If $L/K$ is a finite cyclic extension of degree $n$ with $\Gal(L/K)=\GEN{\sigma}$ and
$a\in K^*$ then $(L/K,\sigma,a)$ denotes the cyclic algebra $L[u|uxu\inv =\sigma(x),
u^n=a]$. Sometimes we abbreviate $(L/K,\sigma,a)$ by writing $(L/K,a)$, if the generator
$\sigma$ is either clear from the context or not relevant. Recall that $(L/K,a)$ is split
if and only if $a$ is a norm of the extension $L/K$ \cite[Theorem 30.4]{Rei}.

A {\it cyclic cyclotomic algebra} is a cyclic algebra $(L/K,a)$ with $L/K$ is a
cyclotomic extension and $a$ is a root of unity. A cyclic cyclotomic algebra $(L/K,a)$ is
a Schur algebra because it is generated over $K$ by the finite metacyclic group
$\GEN{u,\zeta}$, where $\zeta$ is a root of unity of $L$ such that $L=K(\zeta)$.



Quaternion algebras are cyclic algebras of degree $2$ and take the form
$\quat{a,b}{K}=K[i,j|i^2=a,j^2=b,ji=-ij]$, for $a,b\in K^*$. We abbreviate
$\HQ(K)=\quat{-1,-1}{K}$. If $A=\quat{a,b}{K}$ and $\sigma$ is a real embedding of $K$
then $A$ is said to {\it ramify} at $\sigma$ if $\R\otimes_{\sigma(K)} A\simeq \HQ(\R)$,
or equivalently, if $\sigma(a),\sigma(b)<0$. Recall that a totally definite quaternion
algebra is a quaternion algebra $A$ over a totally real field which is ramified at every
real embedding.

Let $G$ be a finite group. A {\it strong Shoda pair} of $G$ is a pair $(M,L)$ of
subgroups of $G$ such that $L\unlhd M \unlhd G$, $M/L$ is cyclic and $M/L$ is maximal
abelian in $N_G(L)/L$. (The definition in \cite{ORS} is more general but for our purposes
we do not need such a generality.) If $M=L$ (and hence $M=G$), then let
$\varepsilon(M,M)=\gorro{M}=\frac{1}{|M|}\sum_{m\in M} m\in \Q M$; otherwise, let
$\varepsilon(M,L) = \prod(\gorro{L}-\gorro{S})$, where $S$ runs on the minimal subgroups
of $M$ containing $L$ properly. Finally, let $e(G,M,L)$ denote the sum of the different
$G$-conjugates of $\varepsilon(M,L)$.

For every positive integer $n$, $\zeta_n$ denotes a complex primitive $n$-th root of
unity. We quote the following result from \cite{ORS}.

\begin{proposition}\label{SSP}
Let $G$ be a finite group and $(M,L)$ a strong Shoda pair of $G$. Let $N=N_{G}(L)$,
$k=[M:L]$ and $n=[G:N]$. Then $e= e(G,M,L)$ is a primitive central idempotent of $\Q G$
and $\Q Ge$ is isomorphic with $M_{n}(\Q(\zeta_k)*^{\sigma}_{\tau} N/M)$, an $n\times
n$-matrix algebra over a crossed product of $N/M$ over the cyclotomic field
$\Q(\zeta_{k})$, with defining action and twisting given as follows: Let $x$ be a
generator of $M/L$ and let $\gamma:N/M\rightarrow N/L$ be a left inverse of the natural
epimorphism $N/L\rightarrow N/M$. Then for every $a,b\in N/M$ one has
    $$ \zeta_k^{\sigma(a)} =
    \zeta_k^i,  \mbox{ if }\, x^{\gamma(a)}= x^i \quad \text{ and } \quad
    \tau(a,b) = \zeta_k^j,  \mbox{ if }\,
    \gamma(ab)\inv \gamma(a)\gamma(b) = x^j.$$
\end{proposition}

In \cite{ORS} it was proved that if $G$ is a finite metabelian groups, then every
primitive central idempotent of $\Q G$ is of the form $e(G,M,L)$ for some strong Shoda
pair $(M,L)$ of $G$. This can be used to compute the Wedderburn decomposition of $\Q G$
for $G$ metabelian. A method to compute the Wedderburn decomposition of $\Q G$ for $G$ an
arbitrary finite group is given in \cite{Olt}. This method has been implemented in the
GAP package \verb+wedderga+ \cite{Wed}. We will make use several times of this method at
different stages of the paper.

We quote the following proposition from \cite{JPRRZ}.

\begin{proposition}\label{ClasKT}
The following are equivalent for a central simple algebra $A$ over a number field $K$.
\begin{enumerate}
\item $A$ is of Kleinian type.
\item $A$ is either a number field or a quaternion algebra which is not ramified at at most
one infinite place.
\item One of the following conditions holds:
\begin{enumerate}
\item $A=K$.
\item $A$ is a totally definite quaternion algebra.
\item $A\simeq M_2(\Q)$.
\item $A\simeq M_2(\Q(\sqrt{d}))$, for $d$ a square-free negative integer.
\item $A$ is a division algebra, $K$ is totally real and $A$ ramifies at all but one real
embeddings of $K$.
\item $A$ is a division algebra, $K$ has exactly one pair of complex (non-real) embeddings
and $A$ ramifies at all real embeddings of $K$.
\end{enumerate}
\end{enumerate}
\end{proposition}

We need the following lemmas.

\begin{lemma}\label{LQuat}
If $K=\Q(\sqrt{d})$ with $d$ a square-free negative integer then
\begin{enumerate}
\item $\HQ(K)$ is a division algebra if and only if $d\equiv 1 \mod 8$.
\item $\quat{-1,-3}{K}$ is a division algebra if and only if $d\equiv 1 \mod 3$.
\end{enumerate}
\end{lemma}

\begin{proof}
(1) Writing $\HQ(K)$ as $(K(\zeta_4)/K,-1)$ one has that $\HQ(K)$ is a division algebra
if and only if $-1$ is a sum of two squares in $K$. It is well known that this is
equivalent to $d\equiv 1 \mod 8$ \cite{FGS}.

(2) Assume first that $A=\quat{-1,-3}{K}$ is not split. Then $A$ is ramified at at least
two finite places $p_1$ and $p_2$. Writing $A$ as $(K(\zeta_3)/K,-1)$ and using
\cite[Theorem~14.1]{Rei}, one deduces that $p_1$ and $p_2$ are divisors of $3$. Thus $3$
is totally ramified in $K$ and this implies that $\left(\frac{D}{3}\right)=1$, where $D$
is the discriminant of $K$ \cite[Theorem~3.8.1]{BS}. Since $D=d$ or $D=4d$ and
$\left(\frac{p}{3}\right)\equiv p \mod 3$, for each rational prime $p$, one has $d =
\left(\frac{d}{3}\right) = \left(\frac{D}{3}\right)\equiv 1 \mod 3$.

Conversely, assume that $d\equiv 1 \mod 3$. Then $3$ is totally ramified in $K$. Let $p$
be a prime divisor of $3$ in $K$. Then the residue field of $K_p$ has order $3$ and
$K_p(\zeta_4)/K_p$ is the unique unramified extension of degree $2$ of $K_p$
\cite[Theorem~5.8]{Rei}. Since $v_p(-3)=1$, we deduce from \cite[Theorem~14.1]{Rei} that
$-3$ is not a norm of the extension $K_p(\zeta_4)/K_p$. Thus $K_p\otimes_K
A=(K_p(\zeta_4)/K_p,-3)$ is a division algebra, hence so is $A$.
\end{proof}

\begin{lemma}\label{DivMeta}
Let $D$ be a division quaternion Schur algebra over a number field $K$. Then $D$ is
generated over $K$ by a metabelian subgroup of $D^*$.
\end{lemma}

\begin{proof}
By means of contradiction we assume that $D$ is not generated over $K$ by a metabelian group. Using Amitsur's
classification of the finite subgroups of division rings (see \cite{A} or \cite{SW}) we deduce that $D$ is generated by
a group $G$ which is isomorphic to one of the following three groups: ${\mathcal O}^*$, the binary octahedral group of
order $48$; $\SL(2,5)$, the binary icosahedral group of order $120$; or $\SL(2,3)\times M$, where $M$ is a metacyclic
group. Recall that ${\mathcal
O}^*=\GEN{x,y,a,b|x^4=x^2y^2=x^2b^2=a^3=1,a^b=a^{-1},x^y=x\inv,x^b=y,x^a=x^{-1}y,y^a=x^{-1}}$.

We may assume without loss of generality that $G$ is one of these three groups. Let $D_1$
be the rational subalgebra of $D$ generated by $G$. It is enough to show that $D_1$ is
generated over $\Q$ by a metabelian group. So we may assume that $D=\Q(G)$ and so $D$ is
one of the factors of the Wedderburn decomposition of $\Q G$.

Computing the Wedderburn decomposition of $\Q {\mathcal O}^*$ and $\Q \, \SL(2,5)$ and
having in mind that $D$ has degree $2$ we obtain that $D\simeq
(\Q(\zeta_8)/\Q(\sqrt{2}),-1)$, if $G={\mathcal O}^*$ and $D\simeq
(\Q(\zeta_5)/\Q(\sqrt{5}),-1)$, if $G=\SL(2,5)$. In both cases $D$ is generated over its
center by a finite metacyclic group.

Finally, assume that $G=\SL(2,3)\times M$, with $M$ metacyclic. Then $D$ is a simple
factor of $A_1\otimes_{\Q} A_2$, where $A_1$ is a simple epimorphic image of $\Q \,
\SL(2,3)$ and $A_2$ is a simple epimorphic image of $\Q M$. If $A_1$ is commutative, then
$D$ is obviously generated by a metabelian group. Assume otherwise that $A_1$ is not
commutative. Computing the Wedderburn decomposition of $\Q \, \SL(2,3)$ we have that
$A_1$ is isomorphic to either $\HQ(\Q)$ or $M_2(\Q(\sqrt{3}))$. Since $D$ is a division
algebra then the second option does not hold. Since $\HQ(\Q)$ is generated over $\Q$ by
the quaternion group of order 8, $A_1\otimes_{\Q} A_2$ is generated over $\Q$ by a direct
product of two metacyclic groups, hence so is $D$.
\end{proof}

For a positive integer $n$ we set $$\eta_n=\zeta_n+\zeta_n^{-1} \quad \mbox{and} \quad
\lambda_n=\zeta_n-\zeta_n^{-1}.$$ Observe that $\eta_n^2-\lambda_n^2=4$ and hence
$\Q(\eta_n^2)=\Q(\lambda_n^2)$. Furthermore, if $n\ge 3$, then $\lambda_n^2$ is totally
negative because if the integer $i$ is relatively prime with $n$ then
$\zeta_n^{2i}+\zeta_n^{-2i}\le -1$. Therefore, if $\lambda_n^2\in K$ then
$\quat{\lambda_n^2,-1}{K}$ ramifies at every real embedding of $K$.

We are ready to classify the Schur algebras of Kleinian type.
\begin{theorem}\label{Clasificacion}
Let $K$ be a number field and let $A$ be a non-commutative central simple
$K$-algebra. Then $A$ is a Schur algebra of Kleinian type if and only if
$A$ is isomorphic to one of the following algebras:
\begin{enumerate}
\item $M_2(K)$, for $K=\Q$ or $\Q(\sqrt{d})$ for $d$ a square-free
negative integer.
\item $\HQ(\Q(\sqrt{d}))$, for $d$ a square-free negative integer, such
that $d\equiv 1 \mod 8$.
\item $\quat{-1,-3}{\Q(\sqrt{d})}$, for $d$ a square-free negative integer, such
that $d\equiv 1 \mod 3$.
\item $\quat{\lambda_n^2,-1}{K}$, where $n\ge 3$,
$\eta_n\in K$ and $K$ has at least one real embedding and at most one pair of complex
(non-real) embeddings.
\end{enumerate}
\end{theorem}

\begin{proof}
That the algebras listed are of Kleinian type follows at once from
Proposition~\ref{ClasKT}.

Now we prove that the algebras listed are cyclic cyclotomic algebras and so they are
Schur algebras. Indeed, if $A=M_2(K)$ as in (1), then $A\simeq (K(\zeta_n)/K,1)$ with
$n=3$ if $K=\Q(\zeta_4)$ and $n=4$ otherwise. If $A$ is a division algebra and either
$A=\HQ(K)$ or $A=\quat{-1,-3}{K}$, for $K=\Q(\sqrt{d})$, then $K$ does not contain
$\zeta_4$ nor $\zeta_3$ (see Lemma~\ref{LQuat}) and hence $\HQ(K)=(K(\zeta_4)/K,-1)$ and
$\quat{-1,-3}{K}=(K(\zeta_3)/K,-1)$. Finally, assume that $A=\quat{\lambda_n^2,-1}{K}$.
Let $B=(K(\zeta_n)/K,-1)=K(\zeta_n)[u|u\zeta_n u\inv = \zeta_n^{-1}, u^2=-1]$, a cyclic
cyclotomic algebra. Then $B=K[\lambda_n,u]$, and $u\lambda_n = - \lambda_n u$. Thus,
$B\simeq \quat{\lambda_n^2,-1}{K}=A$.

Now we prove that if $A$ is a Schur algebra of Kleinian type then one of the cases
(1)-(4) holds. If $A$ is not a division algebra, Proposition~\ref{ClasKT} implies that
$A=M_2(K)$ for $K=\Q$ or an imaginary quadratic extension of $\Q$, so (1) holds.

In the remainder of the proof we assume that $A$ is a division Schur algebra of Kleinian
type. By Lemma~\ref{DivMeta}, $A$ is generated over $K$ by a finite metabelian group $G$.
Then $A=K\otimes_L B$, where $B$ is a simple epimorphic image of $\Q G$ with center $L$
and, by Proposition~\ref{SSP}, $B$ is a cyclic cyclotomic algebra
$(\Q(\zeta_n)/L,\zeta_n^a)$ of degree $2$. Since $A$ is of Kleinian type, so is $B$.

Now we prove that $L$ is totally real. Otherwise, since $L$ is a Galois extension of
$\Q$, $L$ is totally complex and therefore $K$ is also totally complex. By
Proposition~\ref{ClasKT}, both $L$ and $K$ are imaginary quadratic extensions of $\Q$ and
so $L=K$ and $\varphi(n)=4$, where $\varphi$ is the Euler function. Then either (a) $n=8$
and $K=\Q(\zeta_4)$ or $K=\Q(\sqrt{-2})$; or (b) $n=12$ and $K=\Q(\zeta_4)$ or
$K=\Q(\zeta_3)$. If $n=8$, then $B$ is generated over $\Q$ by a group of order $16$
containing an element of order $8$. Since $B$ is a division algebra, $G=Q_{16}$ and so
$B=\HQ(\Q(\sqrt{2}))$, a contradiction. Thus $n=12$ and hence
$B=(\Q(\zeta_{12})/\Q(\zeta_d),\zeta_d^a)$, where $d=6$ or $4$. Since $\zeta_6$ is a norm
of the extension $\Q(\zeta_{12})/\Q(\zeta_3)$, necessarily $d=4$. So
$A=B=(\Q(\zeta_{12})/\Q(\zeta_4),\zeta_4^a)=\quat{\zeta_4^a,-3}{\Q(\zeta_4)}$. Since
$X=1+\zeta_4$, $Y=\zeta_4$ is a solution of the equation $\zeta_4X^2-3Y^2=1$, $\zeta_4$
is a norm of the extension $\Q(\zeta_{12})/\Q(\zeta_4)$, and hence so is $\zeta_4^a$,
yielding a contradiction.

So $L$ is a totally real field of index $2$ in $\Q(\zeta_n)$. Then $L=\Q(\eta_n)$ and
necessarily $B$ is isomorphic to $(\Q(\zeta_n)/\Q(\eta_n),-1)\simeq
\quat{\lambda_n^2,-1}{L}$. This imply that $A\simeq \quat{\lambda_n^2,-1}{K}$. If $K$ has
some embedding in $\R$, then (4) holds. Otherwise $K=\Q(\sqrt{d})$, for some square-free
negative integer $d$. This implies that $L=\Q$. Then $n=3$, $4$ or $6$ and so $A$ is
isomorphic to either $\HQ(K)$ or $\quat{-1,-3}{K}$. Since $A$ is a division algebra,
Lemma~\ref{LQuat} implies that, in the first case, $d\equiv 1 \mod 8$ and condition (2)
holds, and, in the second case, $d\equiv 1 \mod 3$ and condition (3) holds.
\end{proof}

\section{Group algebras}

In this section we classify the group algebras of Kleinian type $KG$ for $K$ a number
field and $G$ a finite group. The classification for $K=\Q$ was given in \cite{JPRRZ}.

We start with some notation. The cyclic group of order $n$ is usually denoted by $C_n$.
To emphasize that $a\in C_n$ is a generator of the group, we write $C_n$ either as
$\GEN{a}$ or $\GEN{a}_n$. Recall that a group $G$ is {\it metabelian} if $G$ has an
abelian normal subgroup $N$ such that $A=G/N$ is abelian. We simply denote this
information as $G=N:A$. To give a concrete presentation of $G$ we will write $N$ and $A$
as direct products of cyclic groups and give the necessary extra information on the
relations between these generators. By $\overline{x}$ we denote the coset $xN$. For
example, the dihedral group of order $2n$ and the quaternion group of order $4n$ can be
given by
    $$\matriz{{llll}
    D_{2n} & = & \GEN{a}_n : \GEN{\overline{b}}_2, & b^2 = 1, a^b = a\inv. \\
    Q_{4n} & = & \GEN{a}_{2n} : \GEN{\overline{b}}_2, & a^b = a\inv, \; b^2 = a^n.}$$
If $N$ has a complement in $G$ then $A$ can be identify with this complement and we write
$G=N\rtimes A$. For example, the dihedral group also can be given by $D_{2n} = \GEN{a}_n
\rtimes \GEN{b}_2$ with $a^b = a\inv$ and the semidihedral groups of order $2^{n+2}$ can
be described as
    $$\matriz{{llll}
    D_{2^{n+2}}^{+} &=& \GEN{a}_{2^{n+1}} \rtimes \GEN{b}_2, &  a^b = a^{2^n + 1}.  \\
    D_{2^{n+2}}^{-} &=& \GEN{a}_{2^{n+1}} \rtimes \GEN{b}_2,  &  a^b = a^{2^n - 1}.}$$

Following the notation in \cite{JPRRZ}, for a finite group $G$, we denote by $\CC(G)$ the
set of isomorphism classes of noncommutative simple quotients of $\Q G$. We generalize
this notation and, for a semisimple group algebra $KG$, we denote by $\CC(KG)$ the set of
isomorphism classes of noncommutative simple quotients of $KG$. For simplicity, we
represent $\CC(G)$ by listing a set of representatives of its elements. For example,
using the isomorphisms
    $$\Q D_{16}^- \cong 4\Q \oplus M_2(\Q) \oplus M_2(\Q(\sqrt{-2})) \quad \mbox{ and } \quad
    \Q D_{16}^+ \cong 4\Q\oplus 2\Q(i) \oplus M_2(\Q(i))$$
one deduces that $\CC(D_{16}^+) = \{M_2(\Q(i))\}$ and $\CC(D_{16}^-) =
\{M_2(\Q),M_2(\Q(\sqrt{-2}))\}$. Notice that $\CC(KG) = \{KZ(A)\otimes_{Z(A)} A: A\in
\CC(G)\}$. If $K$ is a number field then $KG$ is of Kleinian type if and only if each
element of $\CC(KG)$ satisfies one of the conditions of Theorem~\ref{Clasificacion}.

The following groups play an important role in the classification of groups of Kleinian
type.
\medskip

\begin{tabular}{lll}
$\W$ & $=$ & $ \left(\GEN{t}_2 \times \GEN{x^2}_2 \times \GEN{y^2}_2\right) : \left(\GEN{\overline{x}}_2 \times
    \GEN{\overline{y}}_2\right)$, with $t = \left(y,x\right)$ and $Z(\W)=\GEN{x^2,y^2,t}$. \\
$\W_{1n} $ & $=$ & $ \left(\prod\limits_{i=1}^n \GEN{t_i}_2 \times \prod\limits_{i=1}^n \GEN{y_i}_2 \right) \rtimes
    \GEN{x}_4$, with $t_{i} = (y_i,x)$ and $Z(\W_{1n})=\GEN{t_1,\dots,t_n,x^2}$. \\
$\W_{2n} $ & $=$ & $ \left(\prod\limits_{i=1}^n \GEN{y_i}_4 \right) \rtimes \GEN{x}_4$, with $t_{i} = (y_i,x) = y_i^2$
    and $Z(\W_{2n})=\GEN{t_1,\dots,t_n,x^2}$. \\
$\V $ & $=$ & $ \left(\GEN{t}_2 \times \GEN{x^2}_4 \times \GEN{y^2}_4\right) : \left(\GEN{\overline{x}}_2 \times
    \GEN{\overline{y}}_2\right)$, with $t = (y,x)$ and $Z(\W)=\GEN{x^2,y^2,t}$. \\
$\V_{1n} $ & $=$ & $ \left(\prod\limits_{i=1}^n \GEN{t_i}_2 \times \prod\limits_{i=1}^n \GEN{y_i}_4 \right) \rtimes
    \GEN{x}_8$, with $t_{i} = (y_i,x)$ and $Z(\V_{1n})=\GEN{t_1,\dots,t_n,y_1^2,\dots,y_n^2,x^2}$. \\
$\V_{2n} $ & $=$ & $ \left(\prod\limits_{i=1}^n \GEN{y_i}_8 \right) \rtimes \GEN{x}_8$, with $t_{i} = (y_i,x) = y_i^4$
    and $Z(\V_{2n})=\GEN{t_i,x^2}$. \\
$\U_1 $ & $=$ & $ \left(\prod\limits_{1\le i< j \le 3} \GEN{t_{ij}}_2 \times \prod\limits_{ k=1}^3
    \GEN{y_k^2}_2\right): \left(\prod\limits_{k=1}^3 \GEN{\overline{y_k}}_2\right)$, with $t_{ij} = (y_j,y_i)$ and \\
    & & $Z(\U_1)=\GEN{t_{12},t_{13},t_{23},y_1^2,y_2^2,y_3^2}$. \\
$\U_2 $ & $=$ & $ \left(\GEN{t_{23}}_2 \times \GEN{y_1^2}_2 \times \GEN{y_2^2}_4 \times \GEN{y_3^2}_4 \right) :
    \left(\prod\limits_{k=1}^3 \GEN{\overline{y_k}}_2 \right)$, with $t_{ij} = (y_j,y_i)$, $y_2^4=t_{12}$, $y_3^4=t_{13}$
    and \\
    & & $Z(\U_2)=\GEN{t_{12},t_{13},t_{23},y_1^2,y_2^2,y_3^2}$. \\
$\T$ & $=$ & $ \left(\GEN{t}_4 \times \GEN{y}_8\right) : \GEN{\overline{x}}_2$, with $t = (y,x)$ and $x^2 = t^2 = (x,t)$.
    \\
$\T_{1n} $ & $=$ & $ \left(\prod\limits_{i=1}^n \GEN{t_i}_4 \times \prod\limits_{i=1}^n \GEN{y_i}_4 \right) \rtimes
    \GEN{x}_8$, with $t_{i} = (y_i,x)$, $(t_i,x)=t_i^2$ and $Z(\T_{1n})=\GEN{t_1^2,\ldots,t_n^2,x^2}$. \\
$\T_{2n} $ & $=$ & $ \left(\prod\limits_{i=1}^n \GEN{y_i}_8 \right) \rtimes \GEN{x}_4$, with $t_{i} = (y_i,x) =
    y_i^{-2}$ and $Z(\T_{2n})=\GEN{t_1^2,\ldots,t_n^2,x^2}$. \\
$\T_{3n} $ & $=$ & $ \left( \GEN{y_1^2t_1}_2 \times \GEN{y_1}_8\times\prod\limits_{i=2}^n \GEN{y_i}_4 \right) :
    \GEN{\overline{x}}_2$, with $t_i=(y_i,x)$, $(t_i,x)=t_i^2$, $x^2 = t_1^2$, \\
    & & $Z(\T_{3n})=\GEN{t_1^2,y_2^2,\ldots,y_n^2,x^2}$ and, if $i\ge 2$ then
    $t_i=y_i^2$,\\
$\s_{n,P,Q}$& $=$ & $C_3^n\rtimes P = (C_3^n\times Q):\GEN{\overline{x}}_2$, with $Q$ a
subgroup of index $2$ in $P$ and $z^x=z\inv$ for each $z\in C_3^n$.
\end{tabular}
\medskip

We collect the following lemmas from \cite{JPRRZ}.

\begin{lemma}\label{SuperNormal}
Let $G$ be a finite group and $A$ an abelian subgroup of $G$ such that
every subgroup of $A$ is normal in $G$. Let $\mathcal{H}=\{H\mid H \mbox{
is a subgroup of } A \mbox{ with } A/H \mbox{ cyclic and } G'\not\subseteq
H\}$. Then $\CC(G)=\cup_{H\in\mathcal{H}} \CC(G/H)$.
\end{lemma}

\begin{lemma}\label{TimesAbelian}
Let $A$ be a finite abelian group of exponent $d$ and $G$ an arbitrary
group.
\begin{enumerate}
\item If $d|2$ then $\CC(A\times G)=\CC(G)$.
\item If $d|4$ and $\CC(G) \subseteq \left\{M_2(\Q),\HQ(\Q),
\quat{-1,-3}{\Q}, M_2(\Q(\zeta_4))\right\}$ then $\CC(A\times G)\subseteq
\CC(G) \cup \{M_2(\Q(\zeta_4))\}$.
\item If $d|6$ and $\CC(G) \subseteq \left\{M_2(\Q),\HQ(\Q),
\quat{-1,-3}{\Q}, M_2(\Q(\zeta_3))\right\}$ then $\CC(A\times G)\subseteq \CC(G) \cup
\{M_2(\Q(\zeta_3))\}$.
\end{enumerate}
\end{lemma}

\begin{lemma} \label{CCBas}
\begin{enumerate}
\item $\CC(\W_{1n}) = \{M_2(\Q)\}$.
\item $\CC(\W) = \CC(\W_{2n}) = \{M_2(\Q), \HQ(\Q)\}$.
\item $\CC(\V), \CC(\V_{1n}), \CC(\V_{2n}), \CC(\U_1), \CC(\U_2),
\CC(\T_{1n}) \subseteq \{M_2(\Q), \HQ(\Q), M_2(\Q(\zeta_4))\}$.
\item $\CC(\T), \CC(\T_{2n}), \CC(\T_{3n}) \subseteq \{M_2(\Q),
\HQ(\Q), M_2(\Q(\zeta_4)), \HQ(\Q(\sqrt{2})), M_2(\Q(\sqrt{-2}))\}$.
\item Let $G=\s_{n,P,Q}$.
\begin{enumerate}
\item
If $P=\GEN{x}$ is cyclic of order $2^n$ then $\CC(G)=\CC(G/\GEN{x^2})\cup
\left\{\quat{\zeta_{2^{n-1}},-3}{\Q(\zeta_{2^{n-1}})}\right\}$. In particular, if $P=C_2$
then $\CC(G)=\left\{M_2(\Q)\right\}$, if $P=C_4$ then $\CC(G)=\left\{M_2(\Q),
\quat{-1,-3}{\Q}\right\}$ and if $P=C_8$ then $\CC(G)=\left\{M_2(\Q),
\quat{-1,-3}{\Q},M_2(\Q(\zeta_4))\right\}$ .
\item
If $P=\W_{1n}$ and $Q=\GEN{y_1,\ldots,y_n,t_{1},\ldots , t_{n},x^{2}}$ then
$\CC(G)=\left\{M_2(\Q), \quat{-1,-3}{\Q}, M_2(\Q(\zeta_{3}))\right\}$.
\item
If $P=\W_{21}$ and $Q=\GEN{y_{1}^{2},x}$ then $\CC(G)=\{M_2(\Q), \HQ(\Q(\sqrt{3})),
M_2(\Q(\zeta_4)), M_2(\Q(\zeta_{3}))\}$.
\end{enumerate}
\end{enumerate}
\end{lemma}

We are ready to present our classification of the group algebras of Kleinian type.

\begin{theorem}\label{GAKT}
Let $K$ be a number field and $G$ a finite group. Then $KG$ is of Kleinian
type if and only if $G$ is either abelian or an epimorphic image of
$A\times H$ for $A$ an abelian group and one of the following conditions
holds:
\begin{enumerate}
\item $K=\Q$ and one of the following conditions holds.
    \begin{enumerate}
    \item $A$ has exponent $6$ and $H$ is either $\W$, $\W_{1n}$ or $\W_{2n}$, for some
    $n$, or $H=\s_{m,\W_{1n},Q}$ with $Q=\GEN{y_1,\dots,y_m,t_1,\dots,t_m,x^2}$, for some $n$ and $m$.
    \item $A$ has exponent $4$ and $H$ is either $\U_1$, $\U_2$, $\V$, $\V_{1n}$, $\V_{2n}$ or
    $\s_{n,C_8,C_4}$, for some $n$.
    \item $A$ has exponent $2$ and $H$ is either $\T$, $\T_{1n}$, $\T_{2n}$, $\T_{3n}$ or
    $\s_{n,\W_{21},Q}$ with $Q=\GEN{y_1^2,x}$, for some $n$.
    \end{enumerate}

\item $K\neq \Q$ and has at most one pair of complex (non-real) embeddings,
      $A$ has exponent $2$ and $H=Q_8$.

\item $K$ is an imaginary quadratic extension of $\Q$, $A$ has exponent $2$ and $H$ is either
      $\W$, $\W_{1n}$, $\W_{2n}$ or $\s_{n,C_4,C_2}$, for some $n$.

\item $K=\Q(\zeta_3)$, $A$ has exponent $6$ and $H$ is either
      $\W$, $\W_{1n}$ or $\W_{2n}$, for some $n$, or $H=\s_{m,\W_{1n},Q}$
      with $Q=\GEN{y_1,\dots,y_m,t_1,\dots,t_m,x^2}$, for some $n$ and $m$.

\item $K=\Q(\zeta_4)$, $A$ has exponent $4$ and $H$ is either $\U_1, \U_2, \V, \V_{1n},
\V_{2n}, \T_{1n}$ or $\s_{n,C_8,C_4}$, for some $n$.

\item $K=\Q(\sqrt{-2})$, $A$ has exponent $2$ and $H$ is either $D_{16}^-$ or $\T_{2n}$, for some $n$.
\end{enumerate}
\end{theorem}

\begin{proof}
To avoid trivialities we assume that $G$ is non-abelian. The main theorem of \cite{JPRRZ}
states that $\Q G$ is of Kleinian type if and only if $G$ is an epimorphic image of
$A\times H$ for $A$ abelian and $A$ and $H$ satisfy one of the conditions (a)-(c) from
(1). So, in the remainder of the proof, we assume that $K\ne \Q$.

First we prove that if one of the conditions (2)-(6) holds, then $KG$ is of Kleinian
type. Clearly, if $KG$ is of Kleinian type and $H$ is an epimorphic image of $G$, then
$KH$ is also of Kleinian type. So we may assume that $G=A\times H$ with $A$, $H$ and $K$
satisfying one of the conditions listed. We compute $\CC(KG)$ in all the cases.

If (2) holds, then $\CC(KG)=\{\HQ(K)\}$.

If (3) holds then, by Lemmas~\ref{TimesAbelian} and ~\ref{CCBas}, one has $\CC(G)
\subseteq \{M_2(\Q), \HQ(\Q), \quat{-1,-3}{\Q}\}$, and so $\CC(KG)\subseteq \{M_2(K),
\HQ(K), \quat{-1,-3}{K}\}$.

Similarly, if (4) holds then $\CC(G)\subseteq \CC(H)\cup \{M_2(\Q(\zeta_3))\} \subseteq
\{M_2(\Q),\HQ(\Q),\quat{-1,-3}{\Q},M_2(\Q(\zeta_3))\}$. Hence
$\CC(KG)=\{M_2(\Q(\zeta_3))\}$, by Lemma~\ref{LQuat}.

Arguing similarly one deduces that if (5) holds then $\CC(KG)=\{M_2(\Q(\zeta_4))\}$.

Finally, assume that (6) holds. If $H=D_{16}^-$ then
$\CC(G)=\{M_2(\Q),M_2(\Q(\sqrt{-2}))\}$ and so $\CC(KG)=\{M_2(\Q(\sqrt{-2}))\}$.
Otherwise, that is if $H=\T_{2n}$, for some $n$, we show that $\CC(KG) =
\{\HQ(\Q),M_2(\Q),M_2(\Q(\sqrt{-2}))\}$. Let $L$ be a proper subgroup of $H'$ such that
$H'/L$ is cyclic. Using that $(y,x)y^2=1$, for each $y\in \GEN{y_1,\dots,y_n}$ one has
that $\T_{2n}/L$ is an epimorphic image of $\T_{21}\times C_2^{n-1}$. Then
Lemmas~\ref{SuperNormal} and \ref{TimesAbelian} imply that $\CC(G)=\CC(\T_{21})$. So we
may assume that $G=\T_{21}$. Now take $B=Z(\T_{21})=\GEN{t^2,x^2}\simeq C_2^2$ and $L$ a
subgroup of $B$ such that $B/L$ is cyclic. If $t^2\in L$, then $\T_{21}/L$ is an
epimorphic image of $\W$ and therefore $\CC(\T_{21}/L)\subseteq \{M_2(\Q),\HQ(\Q)\}$.
Otherwise $L=\GEN{x^2}$ or $L=\GEN{x^2t^2}$; hence $\T_{21}/L\simeq D_{16}^-$ and so
$\CC(\T_{21}/L)\subseteq \{M_2(\Q),M_2(\Q(\sqrt{-2}))\}$. Using Lemma~\ref{SuperNormal}
one deduces that $\CC(\T_{2n})=\{\HQ(\Q),M_2(\Q),M_2(\Q(\sqrt{-2}))\}$.

Now it is clear that in all the cases the elements of $\CC(KG)$ are of one of the types
(1)-(4) of Theorem~\ref{Clasificacion} and so $KG$ is of Kleinian type.

For the remainder of the proof we assume that $KG$ is of Kleinian type, $G$ is
non-abelian (and $K\ne \Q$). Then $\Q G$ is of Kleinian type, that is $G$ is an
epimorphic image of $A\times H$ for $A$ and $H$ satisfying one of the conditions (a)-(c)
from (1). Furthermore, $K$ has at most one pair of complex embeddings, by
Theorem~\ref{Clasificacion}. We have to show that $K$ and $G$ satisfy one of the
conditions (2)-(6). We consider several cases.

{\it Case 1. Every element of $\CC(G)$ is a division algebra.}

This implies that $G$ is Hamiltonian and so $G\simeq Q_8 \times E \times F$ with $E$
elementary abelian $2$-group and $F$ abelian of odd order \cite[5.3.7]{Rob}. If $F=1$,
then (2) holds. Otherwise, $\CC(KG)$ contains $\HQ(K(\zeta_n))$, where $n$ is the
exponent of $F$. Therefore $n=3$ and $K=\Q(\zeta_3)$, by Theorem~\ref{Clasificacion}.
Since $Q_8$ is an epimorphic image of $\W$, condition (4) holds.

In the remainder of the proof we assume that $\CC(G)$ contains a non-division algebra
$B$. Then $B=M_2(L)$ for some field $L$ and therefore $M_2(KL)\in \CC(KG)$. Since $K\ne
\Q$, $KL$ is an imaginary quadratic extension of $\Q$ and $L\subseteq K$. Let $E$ be the
center of an element of $\CC(G)$. Then $KE$ is the center of an element of $\CC(KG)$. If
$KE\ne K$, then the two complex embeddings of $K$ extends to more than two complex
embeddings of $KE$, yielding a contradiction. This shows that $K$ contains the center of
each element of $\CC(G)$.

{\it Case 2. The center of each element of $\CC(G)$ is $\Q$.}

Then Lemmas~\ref{TimesAbelian} and \ref{CCBas} imply that $\CC(G) \subseteq \{M_2(\Q),
\HQ(\Q), \quat{-1,-3}{\Q}\}$. Using this and the main theorem of \cite{LR} one has that
$G=A\times H$, where $A$ is an elementary abelian $2$-group and $H$ is an epimorphic
image of $\W$, $\W_{1n}$, $\W_{2n}$ or $\s_{n,C_4,C_2}$, for some $n$. So $G$ satisfies
(3).

{\it Case 3. At least one element of $\CC(G)$ has center different from
$\Q$.}

Then the center of each element of $\CC(G)$ is either $\Q$ or $K$. Using
Lemmas~\ref{TimesAbelian} and \ref{CCBas} one has: If $A$ and $H$ satisfy condition (1.a)
then $K=\Q(\zeta_3)$, hence (4) holds. If either $A$ and $H$ satisfy (1.b) or they
satisfy (1.c) with $H=\T_{1n}$, for some $n$, then $K=\Q(\zeta_4)$ and condition (5)
holds.

Otherwise, $A$ has exponent $2$ and $H$ is either $\T, \T_{2n}, \T_{3n}$, for some $n$,
or $\s_{n,\W_{21},Q}$, for $Q=\GEN{y_1^2,x}$. Since $A$ has exponent $2$, one may assume
that $G=A\times H_1$, for $H_1$ an epimorphic image of $H$ and $H_1$ is not an epimorphic
image of any of the groups considered above. We use the standard bar notation for the
images of $\Q H$ in $\Q H_1$.

Assume first that $H=\T$. Then $(M=\GEN{y,t},L=\GEN{ty^{-2}})$ is a strong Shoda pair of
$\T$ and, by using Proposition~\ref{SSP}, one deduces that if
$e=e(\T,M,L)=\widehat{L}\frac{1-y^4}{2}$, then $\Q G e \simeq \HQ(\Q(\sqrt{2}))$. Since
$\HQ(\Q(\sqrt{2}))$ is not of Kleinian type, we have that $\overline{e}=0$, or
equivalently $\overline{y}^4\in \overline{L}$. Hence either $\overline{y}^4=1$,
$\overline{t^2}=1$ or $\overline{t}=\overline{y}^{-2}$. So $H_1$ is an epimorphic image
of either $\T/\GEN{y^4}$, $\T/\GEN{t^2}$ or $\T/\GEN{ty^2}$. In the first case $H_1$ is
an epimorphic image of $\T_{11}$, in the second case $H_1$ is an epimorphic image of
$\V$, which contradicts the hypothesis that $H_1$ is not an epimorphic image of the
groups considered above. Thus $H_1$ is an epimorphic image of $\T/\GEN{ty^2}\simeq
D_{16}^-$. In fact $H_1=D_{16}^-$, because every proper non-abelian quotient of
$D_{16}^-$ is an epimorphic image of $\W$. Then
$\CC(G)=\CC(D_{16}^-)=\{M_2(\Q),M_2(\Q(\sqrt{-2}))\}$ and so $K=\Q(\sqrt{-2})$. Hence
condition (6) holds.

Assume now that $H=\T_{2n}$. By the first part of the proof we have
$\CC(H)=\{\HQ(\Q),M_2(\Q),M_2(\Q(\sqrt{-2}))\}$. Since we are assuming
that one element of $\CC(G)$ has center different from $\Q$, then
$M_2(\Q(\sqrt{-2}))\in \CC(H_1)$ and so $K=\Q(\sqrt{-2})$. Hence condition
(6) holds.

In the two remaining cases we are going to obtain some contradiction.

Suppose that $H=\T_{3n}$. We may assume that $n$ is the minimal positive integer such
that $G$ is an epimorphic image of $\T_{3n}\times A$ for $A$ an elementary abelian
2-group. This implies that $\GEN{\overline{t_1^2},\overline{t_2},\dots,\overline{t_n}}$
is elementary abelian of order $2^n$ and hence
$\GEN{\overline{y_1}^2,\overline{y_2},\dots,\overline{y_n}}\simeq C_4^n$. Let
$M=\GEN{y_1,y_2,\dots,y_n}$ and $L_1=\GEN{t_1y_1^{-2},y_2,y_3,\dots,y_n})$. Then
$(M,L_1)$ is a strong Shoda pair of $H$. By using Proposition~\ref{SSP}, we obtain that
$\Q H e(H,M,L_1) \simeq \HQ(\Q(\sqrt{2}))$. This implies that $\overline{e(H,M,L_1)}=0$,
or equivalently $\overline{t_1^2}=\overline{y_1^4}\in \overline{L_1}$. Since
$\overline{t_1^2}$ has order $2$ and $\overline{t_1^2}\not\in \GEN{t_2,\dots,t_n}$ one
has $\overline{t_1^2}=\overline{t_1y_1^{-2}t_2^{\alpha_2}\cdots t_n^{\alpha_n}}$ for some
$\alpha_1,\dots,\alpha_n\in \{0,1\}$. Since, by assumption, $H_1$ is not an epimorphic
image of $\T_{2n}$, we have $\alpha_i\ne 0$ for some $i\geq 2$. After changing generators
one may assume that $\alpha_2=1$ and $\alpha_i=0$ for $i\ge 3$. Thus
$\overline{t_1}=\overline{y_1^{-2}t_2}$. Let now $L_2 =
\GEN{t_1y_1^{-2},y_2y_1^{-2},y_3,\dots,y_n}$. Then $(M,L_2)$ is also a strong Shoda pair
of $H$ and $\Q Ge(H,M,L_2)\simeq \HQ(\Q(\sqrt{2}))$. The same argument shows that
$\overline{y_1}^4 = \overline{t_1^2} \in \overline{L_2} =
\GEN{\overline{t_1y_1^{-2}},\overline{y_2y_1^{-2}},\overline{y_3},\dots,\overline{y_n}}$.
This yields a contradiction because  $\overline{t_1y_1^{-2}} = (\overline{y_2y_1^{-2}})^2
\in \GEN{\overline{y_2y_1^{-2}},\overline{y_3},\dots,\overline{y_n}}$ and
$\overline{y_1}^4\not\in
\GEN{\overline{y_2y_1^{-2}},\overline{y_3},\dots,\overline{y_n}}$.

Finally assume that $H=\s_{n,\W_{21},Q}$, with $Q=\GEN{y_1^2,x}$ and set $y=y_1$. Since,
by assumption, $G$ does not satisfy (1.a), one has $\overline{t}=\overline{t_1}\ne 1$.
Moreover, as in the previous case one may assume that $n$ is minimal (for $G$ to be a
quotient of $H\times A$, with $A$ elementary abelian $2$-group). Let $M=\GEN{C_3^n,x,t}$
and $L=\GEN{Z_1,tx^2}$, where $Z_1$ is a maximal subgroup of $Z=C_3^n$. Then $(M,L)$ is a
strong Shoda pair of $H$ and $\Q H e(H,M,L)\simeq \HQ(\Q(\sqrt{3}))$. Thus
$0=\overline{e(H,M,L)}=\overline{\widehat{L}(1-\widehat{z})(1-\widehat{t})}$, where $z\in
Z\setminus Z_1$. Comparing coefficients and using the fact that $\overline{t}\ne
\overline{z}$, for each $z\in Z$, we have $\overline{\widehat{L}(1-\widehat{z})}=0$, that
is $\overline{L}=\overline{Z}$ and this contradicts the minimality of $n$.
\end{proof}

\section{Groups of units}

In this section we study the virtual structure of $RG^*$ for $G$ a finite group and $R$
an order in a number field $K$. More precisely, we characterize the finite groups $G$ and
number fields $K$ for which $RG^*$ is finite, virtually abelian, virtually a direct
product of free groups or virtually a direct product of free-by-free groups. We say that
a group {\it virtually} satisfies a group theoretical condition if it has a subgroup of
finite index satisfying the given condition. Notice that the virtual structure of $RG^*$
does not depends on the order $R$ and, in fact, if $S$ is any order in $KG$, then $S^*$
and $RG^*$ are commensurable (see e.g. \cite[Lemma~4.6]{Seh}). Recall that two subgroups
of a given group are said to be {\it commensurable} if their intersection has finite
index in both. It is easy to show that a group commensurable with a free group
(respectively, a free-by-free, a direct product of free groups, a direct product of
free-by-free groups) it is virtually free (respectively, free-by-free, a direct product
of free groups, a direct product of free-by-free groups).

One important tool is the following lemma.

\begin{lemma}\label{RedComp}
Let $A=\prod_{i=1}^n A_i$ be a finite dimensional rational algebra with $A_i$ simple for
every $i$. Let $S$ be an order in $A$ and $S_i$ an order in $A_i$.
\begin{enumerate}
\item $S^*$ is finite if and only if for each $i$, $A_i$ is either $\Q$, an imaginary quadratic extension of $\Q$ or a
totally definite quaternion algebra over $\Q$.
\item $S^*$ is virtually abelian if and only if for each $i$, $A_i$ is either a number field or a totally definite
quaternion algebra.
\item $S^*$ is virtually a direct product of free groups if and only if for each $i$, $A_i$ is either a number field,
a totally definite quaternion algebra or $M_2(\Q)$.
\item $S^*$ is virtually a direct product of free-by-free groups if and only if for each $i$, $S_i^1$ is virtually
free-by-free.
\end{enumerate}
\end{lemma}

\begin{proof}
We are going to use the following facts:
\begin{itemize}
\item[(a)] $S^*$ is commensurable with $\prod_{i=1}^n S_i^*$ and $S_i^*$ is commensurable with $Z(S_i)^*\times S_i^1$.
(This is because $S$ and $\prod_{i=1}^n S_i$ are both orders in $A$.)
\item[(b)] $S_i^1$ is finite if and only if $A_i$ is either a field or a totally definite quaternion algebra
    (see \cite[Lemma 21.3]{Seh} or \cite{KleP}).
\item[(c)] If $A_i$ is neither a field nor a totally definite quaternion algebra then if $S_i^1$ is
commensurable with a direct product of groups $G_1$ and $G_2$ then either $G_1$ or $G_2$
is finite \cite{KR}.
\item[(d)] $S_i^1$ is infinite and virtually free if and only if $A_i\simeq M_2(\Q)$ (see e.g. \cite[page 233]{KleS}).
\end{itemize}

(1) By (a), $S^*$ is finite if and only if $S_i^*$ is finite for each $i$ if and only if
$Z(S_i)^*$ and $S_i^1$ are finite for each $i$. By the Dirichlet Unit Theorem, if $A_i$
is a number field, then $S_i^*$ is finite if and only if $A_i$ is either $\Q$ or an
imaginary quadratic extension of $\Q$. Using this and (b), one deduces that if $A_i$ is
not a number field then $Z(S_i)^*$ and $S_i^1$ are finite if and only if $A_i$ is a
totally definite quaternion algebra over $\Q$.

(2) Clearly $S^*$ is virtually abelian if and only if $S_i^1$ is virtually abelian for
each $i$. If $A_i$ is either a number field or a totally definitive quaternion then
$S_i^*$ is virtually abelian, by (b). Assume otherwise that $A_i$ is neither commutative
nor a totally definite quaternion algebra. Using (c) and the fact that $S_i^1$ is
finitely generated one deduces that $S_i^1$ is virtually abelian, then it is virtually
cyclic and so $A_i=M_2(\Q)$, by (d). This yields a contradiction because $\SL_2(\Z)$
contains a non-abelian free group.

(3) By (a) and \cite[Lemma 3.1]{JR}, $S^*$ is virtually a direct product of free groups
if and only if so is $S_i^1$ for each $i$. As in the previous proof, if $A_i$ is neither
commutative nor a totally definite quaternion algebra, then $S_i^1$ is virtually a direct
product of free groups if and only if $S_i^1$ is virtually free if and only if $A_i\simeq
M_2(\Q)$.

(4) Is proved in \cite[Theorem~2.1]{JPRRZ}.
\end{proof}

The characterization of when $RG^*$ is finite or virtually abelian are easy
generalizations of known results for integral group rings.

\begin{theorem}\label{Finito}
Let $R$ be an order in a number field $K$ and $G$ a finite group. Then $RG^*$ is finite
if and only if one of the following conditions holds:
\begin{enumerate}
\item $K=\Q$ and $G$ is either abelian of exponent dividing $4$ or $6$, or isomorphic to $Q_8\times A$, for $A$ an elementary abelian $2$-group.
\item $K$ is an imaginary quadratic extension of $\Q$ and $G$ is an elementary abelian $2$-group.
\item $K=\Q(\zeta_3)$ and $G$ is abelian of exponent $3$ or $6$.
\item $K=\Q(\zeta_4)$ and $G$ is abelian of exponent $4$.
\end{enumerate}
\end{theorem}

\begin{proof}
If $K=\Q$, then $R=\Z$ and it is well known that $\Z G^*$ is finite if and only if $G$ is
abelian of exponent dividing $4$ or $6$ or it is isomorphic to $Q_8\times A$, for $A$ an
elementary abelian $2$-group.

If one of the conditions (1)-(4) holds, then $KG$ is a direct product of copies of $\Q$,
imaginary quadratic extensions of $\Q$ and $\HQ(\Q)$. Then $RG^*$ is finite by
Lemma~\ref{RedComp}.

Conversely, assume that $RG^*$ is finite and $K\ne \Q$. Then $\Z G^*$ is finite and
therefore $G$ is either abelian of exponent dividing $4$ or $6$ or isomorphic to
$Q_8\times A$, for $A$ an elementary abelian $2$-group. Moreover, $R^*$ is finite and so
$K=\Q(\sqrt{d})$ for $d$ a square-free negative integer. If the exponent of $G$ is $2$
then (2) holds. If the exponent of $G$ is $4$ then one of the simple components of $KG$
is $\Q(\sqrt{d},\zeta_4)$ and therefore $d=-1$, that is (4) holds. If the exponent of $G$
is $3$ or $6$, then one of the simple components of $KG$ is $\Q(\sqrt{d},\zeta_3)$, hence
$d=-3$, and therefore (3) holds.
\end{proof}

\begin{theorem}\label{Vabeliano}
Let $R$ be an order in a number field $K$ and $G$ a finite group. Then $RG^*$ is
virtually abelian if and only if either $G$ is abelian or $K$ is totally real and
$G\simeq Q_8\times A$, for $A$ an elementary abelian $2$-group.
\end{theorem}

\begin{proof}
As in the proof of Theorem~\ref{Finito}, the sufficient condition is a direct consequence
of Lemma~\ref{RedComp}.

Conversely, assume that $RG^*$ is virtually abelian. Then $\Z G^*$ is virtually abelian
and therefore it does not contain a non-abelian free group. Then $G$ is either abelian or
isomorphic to $G\simeq Q_8\times A$, for $A$ an elementary abelian $2$-group \cite{HP}.
In the latter case one of the simple components of $KG$ is $\HQ(K)$ and hence $K$ is
totally real, by Lemma~\ref{RedComp}.
\end{proof}

\begin{theorem}\label{Vfree}
Let $R$ be an order in a number field $K$ and $G$ a finite group. Then $RG^*$ is
virtually a direct product of free groups if and only if either $G$ is abelian or one of the following conditions holds:
\begin{enumerate}

\item $K=\Q$ and $G\simeq H\times A$, for $A$ an elementary abelian $2$-group and $H$ is either $\W$, $\W_{1n}$, $\W_{1n}/\GEN{x^2}$, $\W_{2n}$, $\W_{2n}/\GEN{x^2}$,
$\W_{2n}/\GEN{x^2t_1}$, $\T_{3n}$ or $H=\s_{n,C_{2t},C_t}$, for some $n$ and $t=1,2$ or
$4$.

\item $K$ is totally real and $G\simeq Q_8\times A$, for $A$ an elementary abelian $2$-group.
\end{enumerate}
\end{theorem}

\begin{proof}
The finite groups $G$ such that $\Z G^*$ is virtually a direct product of free groups
were classified in \cite{JR} and coincides with the groups abelian or satisfying condition (1). So, in the remainder of the proof one may assume that $R\ne \Z$, or equivalently
$K\ne \Q$, and we have to show that $RG^*$ is virtually a direct product of free groups
if and only if either $G$ is abelian or (2) holds.

If either $G$ is abelian or (2) holds then $RG^*$ is virtually abelian, hence $RG^*$ is virtually a direct product of free groups, because it is finitely generated.

Conversely, assume that $RG^*$ is virtually a direct product of free groups and $G$ is
non-abelian. Since $K\ne \Q$, $M_2(K)$ is not a simple quotient of $KG$, hence
Lemma~\ref{RedComp} implies that every simple quotient of $KG$ is either a number field
or a totally definite quaternion algebra. In particular, $G$ is Hamiltonian, that is
$G=Q_8 \times A \times F$, where $A$ is an elementary abelian $2$-group and $F$ is abelian of
odd order. If $n$ is the exponent of $F$ then $\HQ(K(\zeta_n))$ is a simple quotient of
$KG$ and this implies that $n=1$ and $K$ is totally real.
\end{proof}

\begin{theorem}\label{Vfreebfree}
Let $R$ be an order in a number field $K$ and $G$ a finite group. Then $RG^*$ is
virtually a direct product of free-by-free groups if and only if either $G$ is abelian or
one of the following conditions holds:
\begin{enumerate}
\item $G$ is an epimorphic image of $A\times H$ with $A$ abelian and $K$,
    $A$ and $H$ satisfy one of the conditions (1), (4), (5) or (6) of Theorem~\ref{GAKT}.
\item $K$ is totally real and $G\simeq A\times Q_8$, for $A$ an elementary abelian $2$-group.
\item $K=\Q(\sqrt{d})$, for $d$ a square-free negative integer, $SL_2(\Z(\sqrt{d}))$ is virtually free-by-free,
and $G\simeq A\times H$ where $A$ is an elementary abelian $2$-group and one of the following holds:
    \begin{enumerate}
    \item $H$ is either $\W_{1n}$, $\W_{1n}/\GEN{x^2}$, $\W_{2n}/\GEN{x^2}$ or $\s_{n,C_2,1}$, for some $n$.
    \item $H$ is either $\W$, $\W_{2n}$ or $\W_{2n}/\GEN{x^2t_1}$, for some $n$ and $d\not\equiv 1 \mod 8$.
    \item $H=\s_{n,C_4,C_2}$ for some $n$ and $d\not\equiv 1 \mod 3$.
    \end{enumerate}
\end{enumerate}
\end{theorem}

\begin{proof} To avoid trivialities we assume that $G$ is non-abelian.

We first show that if $K$ and $G$ satisfy one of the listed conditions then $RG^*$ is
virtually a direct product of free-by-free groups. By Lemma~\ref{RedComp} this is
equivalent to show that if $B\in \CC(KG)$ and $S$ is an order in $B$ then $S^1$ is
virtually free-by-free. By using Lemmas~\ref{TimesAbelian} and \ref{CCBas}, it is easy to
show that if $K$ and $G$ satisfy one of the conditions (1) or (2), then $B$ is either a
totally definite quaternion algebra or isomorphic to $M_2(K)$ for $K=\Q(\sqrt{d})$, with
$d=0,-1,-2$ or $-3$. In the first case $S^1$ is finite and in the second case $S^1$ is
virtually free-by-free (see Lemma~\ref{RedComp} and \cite[Lemma 3.1]{JPRRZ} or
alternatively \cite{KleP}, \cite[page 137]{MR} and \cite{WZ}). If $K$ and $G$ satisfy
condition (3) then Lemmas~\ref{LQuat}, \ref{TimesAbelian} and \ref{CCBas} imply that
$B=M_2(Q(\sqrt{d}))$. Since $S^1$ and $\SL_2(\Z[\sqrt{d}])$ are commensurable and, by
assumption, the latter is virtually free-by-free, we have that $S^1$ is virtually
free-by-free.

Conversely, assume that $RG^*$ is virtually a direct product of free-by-free groups. Let
$B$ be a simple factor of $KG$ and $S$ an order in $B$. By Lemma~\ref{RedComp}, $S^1$ is
virtually free-by-free and hence the virtual cohomological dimension of $S^1$ is at most
$2$. Then $B$ is of Kleinian type by \cite[Corollary~3.4]{JPRRZ}. This proves that $KG$
is of Kleinian type. Furthermore, $B$ is of one of the types (a)-(f) from
Proposition~\ref{ClasKT}. However, the virtual cohomological dimension of $S^1$ is 0, if
$B$ is of type (a) or (b); 1 if $B$ is of type (c); 2 if it is of type (d) or (e); and 3
if $B$ is of type (f) \cite[Remark 3.5]{JPRRZ}. Thus $B$ is not of type (f). Since every
simple factor of $KG$ contains $K$, either $K$ is totally real or $K$ is an imaginary
quadratic extension of $\Q$ and $KG$ is split.

By Theorem~\ref{GAKT}, $G$ is an epimorphic image of $A\times H$ with $A$ abelian and
$K$, $A$ and $H$ satisfying one of the conditions (1)-(6) of Theorem~\ref{GAKT}. If they
satisfy one of the conditions (1), (4), (5) or (6) of Theorem~\ref{GAKT}, then condition
(1) (of Theorem~\ref{Vfreebfree}) holds. So, we assume that $K$, $A$ and $H$ satisfy
either condition (2) or (3) of Theorem~\ref{GAKT}. Since $A$ is elementary abelian, one
may assume that $G=A\times H_1$ with $H_1$ an epimorphic image of $H$.

Assume first that $K$, $A$ and $H$ satisfy condition (2) of Theorem~\ref{GAKT}. Then one
of the simple quotient of $KG$ is isomorphic to $\HQ(K)$. If $K$ is totally real then
condition (2) holds. Otherwise $K=\Q(\sqrt{d})$ for $d$ a square-free negative integer
and $\HQ(K)$ is split. By Lemma~\ref{LQuat}, $d\not\equiv 1 \mod 8$. Since $Q_8\simeq
\W_{21}/\GEN{x^2t_1}$, condition (3b) holds.

Secondly assume that $K$, $A$ and $H$ satisfy condition (3) of Theorem~\ref{GAKT} and set
$K=\Q(\sqrt{d})$ for $d$ a square-free negative integer. Then $\CC(G)\subseteq
\left\{\HQ(\Q),M_2(\Q),\quat{-1,-3}{\Q}\right\}$, by Lemma~\ref{CCBas}. By the main
theorem of \cite{JR}, $H_1$ is isomorphic to either $\W$, $\W_{1n}$, $\W_{2n}$,
$\W_{1n}/\GEN{x_1^2}$, $\W_{2n}/\GEN{x_1^2}$, $\W_{2n}/\GEN{x_1^2t_1}$, $\s_{n,C_2,1}$ or
$\s_{n,C_4,C_2}$. If $H$ is either $\W_{1n}$, $\W_{1n}/\GEN{x_1^2}$,
$\W_{2n}/\GEN{x_1^2}$ or $\s_{n,C_2,1}$, then (3a) holds. If $H$ is either $\W$,
$\W_{2n}$ or $\W_{2n}/\GEN{x^2t_1}$ then $\CC(G)=\{M_2(\Q),\HQ(\Q)\}$ and, arguing as in
the previous paragraph, one deduces that $d\not\equiv 1 \mod 8$. In this case condition
(3b) holds. Finally, if $G=\s_{n,C_4,C_2}$ then $\CC(G)=\{M_2(\Q),\quat{-1,-3}{\Q}\}$ and
using the second part of Lemma~\ref{LQuat} one deduces that $d\not\equiv 1 \mod 3$,  and
condition (3c) holds.
\end{proof}

The main theorem of \cite{JPRRZ} states that a finite group $G$ is of Kleinian type if
and only if $\Z G^*$ is commensurable with a direct product of free-by-free groups. One
implication is still true when $\Z$ is replace by an arbitrary order in a number field. This is a
consequence of Theorems~\ref{GAKT} and \ref{Vfreebfree}.

\begin{corollary}\label{ffK}
Let $R$ be an order in a number field $K$ and $G$ a finite group. If $RG^*$ is
commensurable with a direct product of free-by-free groups then $KG$ is of Kleinian type.
\end{corollary}

It also follows from Theorems~\ref{GAKT} and \ref{Vfreebfree} that the converse of
Corollary~\ref{ffK} fails. The group algebras $KG$ of Kleinian type for which the group
of units of an order in $KG$ is not virtually a direct product of free-by-free groups
occur under the following circumstances, where $G=A\times H$ for $A$ an elementary
abelian $2$-group:
\begin{enumerate}
\item $K$ is a number field with exactly one pair of complex embeddings and at least one real embedding and $H=Q_8$.
\item $K=\Q(\sqrt{d})$, for $d$ square free-negative integer with $d\equiv 1 \mod 8$ and $H=\W$, $\W_{2n}$ or $\W_{2n}/\GEN{x^2t_1}$, for some $n$.
\item $K=\Q(\sqrt{d})$ for $d$ a square-free negative integer with $d\equiv 1 \mod 3$ and $H=\s_{n,C_4,C_2}$.
\item $K=\Q(\sqrt{d})$ and $d$ and $H$ satisfy one of the conditions (3a)-(3c) from Theorem~\ref{Vfreebfree},
    but $\SL_2(\Z[\sqrt{d}])$ is not virtually free-by-free.
\end{enumerate}

Thus, we have a good description of the virtual structure of $RG^*$ for $R$ an order in
$KG$ if $KG$ is of Kleinian type, except in the four cases above. It has been conjecture
that $\SL_2(\Z[\sqrt{d}])$ is virtually free-by-free for every negative integer. This
conjecture has been verified for $d=-1,-2,-3,-7$ and $-11$. Thus, maybe the last case
does not occur and the hypothesis of $\SL_2(\Z[\sqrt{d}])$ being virtually free-by-free
in Theorem~\ref{Vfreebfree} is superfluous.

Notice that $K=\Q(\sqrt{-7})$ and $G=Q_8=\W_{11}/\GEN{x^2t_1}$ is an instance of case (2)
above and, if $R$ is an order in $K$, then $RQ_8^*$ is commensurable with $\HQ(R)^*$. A
presentation of $\HQ(R)^*$, for $R$ the ring of integers of $\Q(\sqrt{-7})$ has been
computed in \cite{CEGA}.

\bibliographystyle{amsalpha}

\end{document}